\documentclass[12pt]{article}
\usepackage{latexsym}
\usepackage{amssymb}
\newcommand{\ga}{\alpha}

\renewcommand{\gg}{\gamma}
\newcommand{\gd}{\delta}

\newcommand{\gk}{\kappa}
\newcommand{\gl}{\lambda}

\newcommand{\go}{\omega}

\newcommand{\la}{\langle}
\newcommand{\ra}{\rangle}
\newcommand{\ov}{\overline}
\newcommand{\add}{{\rm Add}}
\newcommand{\coll}{{\rm Coll}}

\newcommand{\card}[1]{{\vert #1 \vert} }

\newcommand{\forces}{\Vdash}

\newcommand{\FP}{{\mathbb P}}
\newcommand{\FQ}{{\mathbb Q}}
\newcommand{\FR}{{\mathbb R}}

\newtheorem{theorem}{Theorem}
\newtheorem{corollary}[theorem]{Corollary}

\newenvironment{proof}{\noindent{\bf
Proof:}}{\nopagebreak\mbox{}\newline
\makebox[\textwidth]{\hfill$\square$}
\par\bigskip}

\def\P{{\mathbb P}}

\newcommand{\restrict}{\upharpoonright}

\def\cof{\mathop{\rm cof}\nolimits}

\def\intersect{\cap}

\def\and{\mathrel{\kern1pt\&\kern1pt}}
\def\image{\mathbin{\hbox{\tt\char'42}}}

\def\<#1>{\langle\,#1\,\rangle}

\newcommand{\cp}{\mathop{\rm cp}}
\title{Large Cardinals with Few Measures\thanks{
         cardinal, strongly compact cardinal, measurable cardinal,
         normal measure.}}
\author{{Arthur W.~Apter}\thanks{The research of the first and
        third author was partially supported by PSC-CUNY grants and
        CUNY Collaborative Incentive Grants.}\\
        {Department of Mathematics}\\
        {Baruch College of CUNY}\\
        {New York, New York 10010}\\
        {http://faculty.baruch.cuny.edu/apter}\\
        {awabb@cunyvm.cuny.edu}\\
        \\
        {James Cummings}\thanks{The second author's research was partially
        supported by NSF Grant DMS-0400982.}\\
        {Department of Mathematical Sciences}\\
        {Carnegie Mellon University}\\
        {Pittsburgh, Pennsylvania 15213}\\
        {http://www.math.cmu.edu/users/jcumming}\\
        {jcumming@andrew.cmu.edu}\\
        \\
        Joel David Hamkins${}^{\dag}$\\
        The City University of New York\\
        The College of Staten Island of CUNY and\\
        The CUNY Graduate Center, Mathematics\\
        365 Fifth Avenue, New York, New York 10016\\
        http://jdh.hamkins.org\\
        jdh@hamkins.org}
\begin{document}\maketitle
\begin{abstract}
We show, assuming the consistency of one measurable
cardinal, that it is consistent for there to be exactly
$\gk^+$ many normal measures on the least measurable
cardinal $\gk$. This answers a question of Stewart Baldwin.
The methods generalize to higher cardinals, showing that
the number of $\gl$ strong compactness or $\gl$
supercompactness measures on $P_\gk(\gl)$ can be exactly
$\gl^+$, if $\gl > \gk$ is a regular cardinal.
%the number of $\lambda$ supercompactness or $\lambda$ strong
%compactness measures over $P_\kappa(\lambda)$ can be exactly $\lambda^+$.
We conclude with a list of open questions. Our proofs use a
critical observation due to James Cummings.
\end{abstract}
\baselineskip=24pt
%\newpage
\section{Introduction and Preliminaries}\label{s1}
Set theorists have long been occupied with the problem of
determining the number of normal measures there can be on a
measurable cardinal. In \cite{Ku}, Kunen showed that it is
consistent, relative to the existence of a measurable
cardinal, for there to be exactly one normal measure on the
least measurable cardinal $\gk$. In \cite{KP}, Kunen and
Paris showed that it is consistent, relative to the
existence of a measurable cardinal, for there to be exactly
$2^{2^\gk}$ many normal measures on the least measurable
cardinal $\gk$, the maximal possible number. In the result
of \cite{KP}, it can be the case that $2^{2^\gk} =
\gk^{++}$. In \cite{Mi}, Mitchell showed that it is
consistent, relative to a measurable cardinal $\gk$ with
$o(\gk) = \gd$,
%(for $\gd \le \gk^{++}$ an arbitrary finite or infinite cardinal),
for there to be exactly $\gd$ many normal measures on
$\gk$. In this result, $\gd \le \gk^{++}$ is an arbitrary
finite or infinite cardinal. In \cite{Ba}, S$.$ Baldwin
generalized Mitchell's results of \cite{Mi} and showed that
it is consistent, relative to measurable cardinals of high
Mitchell order, for there to be exactly $\gd$ many normal
measures on the least measurable cardinal $\gk$, where $\gd
< \gk$ is an arbitrary finite or infinite cardinal. Note
that the result of \cite{KP} uses forcing, while the
results of \cite{Ku}, \cite{Mi}, and \cite{Ba} use inner
model techniques, so that the GCH holds in the models
constructed.

Until recently, little had been known concerning
generalizations of these results to strongly compact and
supercompact cardinals, primarily because of the limited
inner model theory available for these cardinals. Some
questions remained open even at the level of measurable
cardinals. For instance, Baldwin \cite{Ba} leaves open the
questions of whether it is consistent, relative to some
large cardinal hypothesis, for there to be either exactly
$\gk$ many or exactly $\gk^+$ many normal measures on the
least measurable cardinal $\gk$.

In this paper, we rectify the situation described in the
preceding paragraph to some extent by proving the following
theorem.

\newtheorem{maintheorem}[theorem]{Main Theorem}
\begin{maintheorem}\label{MainTheorem} If $\kappa$ is a
measurable cardinal, then there is a forcing extension,
neither creating nor destroying any measurable cardinals,
where there are exactly $\kappa^+$ many normal measures on
$\gk$.
\end{maintheorem}

\noindent Consequently, there can be exactly $\kappa^+$
many normal measures on the least measurable cardinal
$\gk$. The argument used in the proof of Theorem
\ref{MainTheorem} easily adapts to show that $\kappa$ can
have a specified Mitchell order, or have a Laver function
for measurability, or be $\mu$-measurable, while still
having exactly $\kappa^+$ many normal measures.

This argument also extends to the case of supercompactness
and strong compactness. For example, Theorem
\ref{Theorem.Supercompact} shows that for $\gl > \gk$ a
regular cardinal, it is relatively consistent that $\gk$ is
$\gl$ supercompact, but there are fewer than the maximal
number of fine, normal $\gk$-additive measures on
$P_\gk(\gl)$. In addition, Theorem
\ref{Theorem.StronglyCompact} shows the same for strong
compactness.

All of our results are proved via forcing. The core
observation in the proof of the Main Theorem was made by
James Cummings and subsequently adapted to the various
large cardinal contexts.

We take this opportunity to mention some preliminary
information. For $\gk < \gl$ cardinals, $\gk$ regular,
${\rm Coll}(\gk, \gl)$ is the standard L\'evy collapse of
$\gl$ to $\gk$. For $\gk$ a regular cardinal and $\gg$ an
ordinal, $\add(\gk, \gg)$ is the standard partial ordering
for adding $\gg$ Cohen subsets of $\gk$.
%${\rm Add}(\go, 1)$ is the
%standard partial ordering
%for adding a single Cohen
%subset of $\go$.
The partial ordering $\FP$ is {\em $\gk$-distributive} if
the intersection of $\gk$ many dense open subsets of $\FP$
is dense open. The partial ordering $\FP$ is {\em
$\gk$-strategically closed} if player II has a winning
strategy, ensuring that play can continue for $\kappa$ many
steps, in the two person game in which the players
construct a decreasing sequence of conditions $\la p_\ga
\mid \ga \le \gk \ra$, where player I plays at all odd
stages and player II plays at all even stages (including
all limit stages), with the trivial condition played at
stage 0. The partial ordering $\FP$ is {\em
$\gk^+$-directed closed} if any directed collection of
conditions of cardinality at most $\gk$ has a common
extension. Such partial orderings are necessarily
$\gk$-strategically closed.

A forcing notion $\FP$ (and the forcing extensions to which
it gives rise) admits a {\em closure point} at $\delta$ if
it factors as $\FQ*\dot \FR$, where $\FQ$ is nontrivial,
$\card{\FQ} \leq\delta$, and $\forces_\FQ ``\dot\FR$ is
$\delta$-strategically closed''. Our arguments will rely on
the following consequence of the main result of \cite{H5}
(which generalizes results of \cite{H2}).

\begin{theorem}[\cite{H5}]\label{Theorem.ClosurePoint}
If $V \subseteq V[G]$ admits a closure point at $\delta$
and $j:V[G]\to M[j(G)]$ is an ultrapower embedding in
$V[G]$ with $\delta<\cp(j)$, then $j\restrict V:V\to M$ is
a definable class in $V$.
\end{theorem}

This theorem follows from \cite[Theorem 3, Corollary
14]{H5}. If $j:V[G]\to M[j(G)]$ witnesses the $\lambda$
supercompactness of $\kappa$ in $V[G]$, then by
\cite[Corollary 4]{H5}, the restriction $j\restrict V:V\to
M$ witnesses the $\lambda$ supercompactness of $\kappa$ in
$V$. In the case of strong compactness, a slight additional
hypothesis is used. Specifically, a forcing extension
$V\subseteq V[G]$ exhibits {\em $\kappa$-covering} if every
set of ordinals $x$ of size less than $\kappa$ in $V[G]$ is
covered by a set of ordinals $y$ of size less than $\kappa$
in $V$ (this captures the power of {\em mildness} in
\cite{H2}). The fact proved in \cite[Theorem 31]{H5} is
that if $V\subseteq V[G]$ has $\kappa$-covering and admits
a closure point at $\delta<\kappa$, and $j:V[G]\to M[j(G)]$
is the ultrapower by a fine, $\gk$-additive measure on
$P_\kappa(\lambda)$ in $V[G]$, then $j\restrict V:V\to M$
is a definable class in $V$ and witnesses the $\lambda$
strong compactness of $\kappa$ in $V$.

Finally, let us mention that we assume familiarity with the
large cardinal properties of measurability, strong
compactness, and supercompactness, along with some other
related notions. Interested readers may consult \cite{J}
for further details.

\section{Limiting the Number of Measures
on a Measurable Cardinal} \label{Section.Measurable}

We now prove the Main Theorem.
%which we state here with a more detailed conclusion.
%
%\begin{theorem}\label{Theorem.Detailed}\label{t1}
%If $\kappa$ is a measurable cardinal, then there is a
%forcing extension $V^\FP$, neither creating nor destroying
%any measurable cardinals, where $\kappa$ carries exactly
%$\kappa^+$ normal measures. Furthermore, if $\kappa$
%has at least $\kappa^+$ normal measures in $V$, then
%every normal ultrapower in $V$ lifts uniquely to the
%extension and every normal ultrapower in the extension
%arises in this way. Consequently, $o(\kappa)$ is preserved.
%\end{theorem}
%

\newcounter{oldtheoremcount}
\setcounter{oldtheoremcount}{\value{theorem}}
\setcounter{theorem}{0}
\begin{maintheorem} If $\kappa$ is a
measurable cardinal, then there is a forcing extension,
neither creating nor destroying any measurable cardinals,
where there are exactly $\kappa^+$ many normal measures on
$\gk$.
\end{maintheorem}
\setcounter{theorem}{\value{oldtheoremcount}}

\begin{proof}
Suppose that $\kappa$ is a measurable cardinal in $V$. By forcing if necessary, we may assume that there are at least $\kappa^+$ many normal measures
on $\gk$. One way to accomplish this, for example, is via a reverse Easton iteration using $\add(\gd^+,1)$ at every inaccessible cardinal
$\gd\leq\kappa$, and trivial forcing at all other stages. Note that this forcing ensures the GCH will hold in the extension at all such nontrivial
stages of forcing $\gg\le\gk$. Standard arguments (see Lemma 1.1 of \cite{A01} or Lemma 6 of \cite{C93}) then show that every measurable cardinal
$\gg\le\gk$ in $V$ remains measurable after the forcing, with $2^{\gg^+}$ many normal measures, computed in the extension. The results of \cite{LS}
show that all measurable cardinals greater than $\gk$ are also preserved. Further, an easy application of Theorem \ref{Theorem.ClosurePoint} shows
that this iteration creates no new measurable cardinals. So we assume without loss of generality that there are at least $\kappa^+$ many normal
measures on $\gk$ in $V$.

Let $\P=\add(\omega,1)\ast \dot
\coll(\kappa^+,2^{2^{\kappa}})$ be the forcing to add a
Cohen real and then collapse $2^{2^\kappa}$ to $\kappa^+$,
and suppose that $V[c][G]$ is the resulting forcing
extension. Every normal measure on $\kappa$ in $V$
generates a unique normal measure in $V[c]$, since this is
small forcing (see \cite{LS}). These remain normal measures
in $V[c][G]$, since no additional subsets of $\kappa$ are
added by the collapse forcing. Thus, there are at least
$\kappa^+$ many normal measures on $\gk$ in $V[c][G]$.
Conversely, suppose that ${\cal U}$ is a normal measure on
$\kappa$ in $V[c][G]$, with the associated ultrapower
embedding $j:V[c][G]\to M[c][j(G)]$. In particular,
$X\in{\cal U}$ if and only if $\kappa\in j(X)$ for all $X
\subseteq \kappa$ in $V[c][G]$. By Theorem
\ref{Theorem.ClosurePoint}, it follows that the restriction
$j\restrict V:V\to M$ is a definable class in $V$. Since
the forcing to add $c$ is small with respect to $\kappa$,
it follows (by \cite{LS}) that $j\restrict V$ lifts
uniquely to $V[c]$, and so $j\restrict V[c]:V[c]\to M[c]$
is a definable class in $V[c]$. The key observation is now
that because $V[c]$ and $V[c][G]$ have the same subsets of
$\kappa$, one can reconstruct ${\cal U}$ inside $V[c]$ by
observing $X\in{\cal U}$ if and only if $\kappa\in j(X)$,
using only $j\restrict V[c]$. Thus, ${\cal U}\in V[c]$. So
every normal measure on $\kappa$ in $V[c][G]$ is actually
in $V[c]$. The number of such normal measures, therefore,
is at most $(2^{2^\kappa})^{V[c]}$, which is $\kappa^+$ in
$V[c][G]$, because $(2^{2^\kappa})^{V[c]}$ was collapsed by
$G$. So in $V[c][G]$, there are exactly $\kappa^+$ many
normal measures on $\gk$, as desired. By the results of
\cite{LS} and the closure properties of the L\'evy
collapse, forcing with $\FP$ doesn't affect measurable
cardinals either above or below $\gk$, so this completes
the proof of Theorem \ref{MainTheorem}.
\end{proof}

The proof of the Main Theorem is easily adapted, as we now
illustrate. In particular, our first corollary answers one
of Baldwin's questions from \cite{Ba}, by showing that it
is consistent for there to be exactly $\gk^+$ many normal
measures on the least measurable cardinal $\gk$.

\begin{corollary}\label{Corollary.LeastMeasurable}
If $\kappa$ is measurable, then there is a forcing
extension in which $\kappa$ is the least measurable
cardinal and there are exactly $\kappa^+$ many normal
measures on $\gk$.
\end{corollary}

\begin{proof}
Suppose that $\kappa$ is measurable. By either iterating
Prikry forcing (see \cite{Ma}) or iterating non-reflecting
stationary set forcing (see \cite{AC1}) if necessary, we
may arrange that $\kappa$ becomes the least measurable
cardinal.
%and has at least $\kappa^+$ normal measures.
By Theorem \ref{MainTheorem}, there is a further extension
in which $\kappa$ remains the least measurable cardinal and
there are exactly $\kappa^+$ many normal measures on
$\gk$.\end{proof}

Our next corollary shows that a measurable cardinal $\gk$
of nontrivial Mitchell rank can have exactly $\gk^+$ many
normal measures.

\begin{corollary}\label{Order}

Suppose $\gk$ is measurable and $o(\gk) = \gd$. There is
then a class model in which $\gk$ is measurable, $o(\gk) =
\gd$, and there are exactly $\gk^+$ many normal measures on
$\gk$.

\end{corollary}

\begin{proof}
By the main result of \cite{C93}, we may assume without
loss of generality (by first passing to a Mitchell inner
model and then doing the relevant forcing if necessary)
that in our ground model $V$, $o(\gk) = \gd$ and there are
exactly $\gk^{++}$ many normal measures on $\gk$.
Alternatively, one can accomplish this purely by forcing:
first increase the number of normal measures on $\kappa$
with a reverse Easton iteration, as in the Main Theorem,
and then observe that because all the ground model normal
ultrapower embeddings lift to the extension, we preserve
$o(\gk)\geq\delta$; if it happens that $o(\gk)>\delta$ in
the extension, then $o(\gk)=\delta$ in the ultrapower by a
normal measure of rank $\delta$, where there are still
sufficient normal measures.

We now force over $V$ with the partial ordering $\FP$ of
Theorem \ref{MainTheorem}. If $j:V\to M$ is any normal
ultrapower embedding for $\kappa$ in $V$, then it lifts
uniquely through the small forcing $\add(\go, 1)$ to
$j:V[c]\to M[c]$, and then uniquely again through the
directed closed forcing $\coll(\kappa^+,2^{2^{\kappa}})$ to
$j:V[c][G]\to M[c][j(G)]$, since $j(G)$ can (and must) be
taken to be the filter generated by $j\image G$.
Conversely, we have already argued above that if
$j:V[c][G]\to M[c][j(G)]$  is the ultrapower by a normal
measure ${\cal U}$ in $V[c][G]$, then $j:V[c]\to M[c]$ is
definable in $V[c]$. In fact, $j\restrict V[c]$ is the
ultrapower by $\cal U$ in $V[c]$, since the collapse
forcing $G$ adds no new functions from $\kappa$ to $V[c]$.
Since $V[c]$ is a small forcing extension, it follows that
$j\restrict V$ is the ultrapower in $V$ by ${\cal
U}\intersect V\in V$. So we have established that every
normal ultrapower embedding in $V$ lifts to $V[c][G]$, and
all ultrapower embeddings in $V[c][G]$ arise in this way.
It follows that $o(\kappa)$ is preserved, since by
induction, the Mitchell rank of every measure is preserved
to its unique extension in $V[c]$.
\end{proof}

By the definition of $\FP$, the cardinal ${(\gk^{++})}^V$
is of course collapsed. Thus, if $o(\gk) = \gk^{++}$ in
$V$, then in our forcing extension, the Mitchell order of
$\gk$ remains $(\kappa^{++})^V$, which is now an ordinal
between $\gk^+$ and $\gk^{++}$ of the extension.

Next, we adapt the argument to allow for Laver functions.
As in \cite{Hamkins:LaverDiamond}, define that
$\ell:\kappa\to V_\kappa$ is a {\em Laver function for
measurability} if for every $x\in H_{\kappa^+}$, there is a
normal ultrapower embedding $j:V\to M$ with $\cp(j)=\kappa$
and $j(\ell)(\kappa)=x$. Since different values for $x$
give rise to different induced normal measures, the
existence of a Laver function for measurability implies
that there are at least $2^\kappa=|H_{\kappa^+}|$ many
normal measures on $\gk$. Hamkins asked in
\cite{Hamkins:LaverDiamond} whether or not the existence of
such a Laver function implies that $\kappa$ must have
$2^{2^\kappa}$ many normal measures. This is answered in
the negative by the following result.

\begin{corollary}\label{Corollary.LaverFunction}
If $\kappa$ is a measurable cardinal, then there is a
forcing extension in which there are exactly $\kappa^+$
many normal measures on $\gk$, yet $\gk$ has a Laver
function for measurability.
\end{corollary}

\begin{proof} Suppose that $\kappa$ is measurable. We may
assume, by preliminary forcing as in \cite[Theorem 2.3]{H4}
if necessary, that $\kappa$ already has in $V$ a Laver
function $\ell:\kappa\to V_\kappa$. Thus, in $V$, there are
at least $\kappa^+$ many normal measures on $\kappa$.

Let $V[c][G]$ be the forcing extension of Theorem
\ref{MainTheorem}, where there are exactly $\kappa^+$ many
normal measures on $\gk$. Define
$\ell^*(\alpha)=\ell(\alpha)_c$, provided $\ell(\alpha)$ is
an $\add(\omega,1)$-name (choose anything otherwise). For
any $x\in H_{\kappa^+}^{V[c][G]}=H_{\kappa^+}^{V[c]}$,
there is a name $\dot x\in H_{\kappa^+}^V$ such that
$x=\dot x_c$. Since $\ell$ is a Laver function in $V$,
there is a normal ultrapower embedding $j:V\to M$ with
$j(\ell)(\kappa)=\dot x$. This embedding lifts (uniquely)
to $j:V[c][G]\to M[c][j(G)]$. To conclude the argument,
observe that $j(\ell^*)(\kappa)=j(\ell)(\kappa)_c=\dot
x_c=x$, so $\ell^*$ is a Laver function in $V[c][G]$.
\end{proof}

A cardinal $\kappa$ is {\em $\mu$-measurable} if there is
an embedding $j:V\to M$ with critical point $\kappa$ such
that the induced normal measure ${\cal
U}=\{X\subseteq\kappa\mid \kappa\in j(X)\}$ is in $M$.

\begin{corollary}\label{Corollary.MuMeasurable}
If $\kappa$ is $\mu$-measurable, then there is a forcing
extension preserving this in which there are exactly
$\kappa^+$ many normal measures on $\gk$.
\end{corollary}

\begin{proof}
Suppose that $\kappa$ is $\mu$-measurable. From this, it
follows that there are at least $\kappa^+$ many normal
measures (of varying Mitchell rank) on $\gk$. Thus, in the
extension $V[c][G]$ of Theorem \ref{MainTheorem}, there are
exactly $\kappa^+$ many normal measures on $\gk$.

The $\mu$-measurability of $\kappa$ in $V$ is witnessed by
an ultrapower embedding $j:V\to M$ by a (non-normal)
measure $\nu$ on $\kappa$, with the induced normal measure
${\cal U}$ in $M$. Just as in Theorem \ref{MainTheorem},
this embedding lifts uniquely to an embedding
$j^*:V[c][G]\to M[c][j(G)]$ in $V[c][G]$. This embedding
witnesses $\mu$-measurability in $V[c][G]$, since the
induced normal measure $\cal U^*$ is precisely the measure
generated as a filter by ${\cal U}\in M$.
\end{proof}

One can easily generalize the Main Theorem to cardinals
other than $\kappa^+$. For example, if $\kappa$ is
measurable, $\delta$ is a regular cardinal in the interval
$(\kappa^+,2^{2^\kappa}]$ and there are at least $\delta$
many normal measures on $\gk$, then in the forcing
extension $V[c][G]$, obtained by forcing with
$\add(\omega,1) \ast \dot \coll(\delta,2^{2^\kappa})$,
there will be exactly $\delta$ many normal measures on
$\kappa$. And as above, one can similarly preserve various
properties of $\kappa$, such as $\mu$-measurability, the
existence of a Laver function, or any particular value of
$o(\kappa)$.

\section{Limiting the Number of Measures
on Higher Cardinals}\label{s3}

We now extend the method to the case of supercompactness
and strong compactness.

\begin{theorem} \label{Theorem.Supercompact}
Suppose that $\kappa$ is $\lambda$ supercompact,
$\lambda>\kappa$ is regular, and the GCH holds. Then there
is a forcing extension preserving this in which there are
exactly $\lambda^+$ many fine, normal, $\gk$-additive
measures on $P_\kappa(\lambda)$.
\end{theorem}

\begin{proof} By forcing if necessary, we may
%first pump up the number of measures
assume without loss of generality (see, for example,
\cite[Theorem1.15]{H4}) that in $V$, there are at least
$\lambda^+$ many fine, normal, $\gk$-additive measures on
$P_\kappa(\lambda)$. Next, we force as in Theorem
\ref{MainTheorem} with $\P=\add(\omega,1) \ast \dot
\coll(\lambda^+,2^{2^\lambda})$, giving rise to the forcing
extension $V[c][G]$. Standard arguments show that the GCH
holds in $V[c][G]$.

As before, every $\lambda$ supercompactness embedding
$j:V\to M$ in $V$ lifts uniquely through the small forcing
to $j:V[c]\to M[c]$, and then uniquely to the full
extension $j:V[c][G]\to M[c][j(G)]$. This is because the
collapse forcing is $\lambda$-distributive, which implies
that the filter generated by $j\image G$ is already
$M[c]$-generic, and so one can (and must) take $j(G)$ to be
this filter. Conversely, if ${\cal U}$ is a fine, normal,
$\gk$-additive measure on $P_\kappa(\lambda)$ in $V[c][G]$,
with the associated ultrapower embedding $j:V[c][G]\to
M[c][j(G)]$, then by Theorem \ref{Theorem.ClosurePoint},
the restricted embedding $j\restrict V:V\to M$ is a
definable class in $V$. This embedding lifts uniquely to
$j\restrict V[c]:V[c]\to M[c]$ in $V[c]$.

The point now is just as in Theorem \ref{MainTheorem},
namely that $V[c]$ and $V[c][G]$ have the same subsets of
$P_\kappa(\lambda)$. Therefore, $\cal U$ is constructible
from $j\restrict V[c]$ in $V[c]$, because $X\in {\cal U}$
if and only if $j\image\lambda\in j(X)$. Hence, every
$\lambda$ supercompactness measure in $V[c][G]$ is actually
in $V[c]$. Consequently, the number of such measures is at
most $(2^{2^\lambda})^{V[c]}$, which has size $\lambda^+$
in $V[c][G]$.
\end{proof}

\begin{theorem} \label{Theorem.StronglyCompact}

Suppose that $\kappa$ is $\lambda$ supercompact, $\lambda >
\kappa$ is regular, and the GCH holds. Then there is a
forcing extension in which $\kappa$ is $\lambda$ strongly
compact but not $\lambda$ supercompact, the GCH holds, and
there are exactly $\gl^+$ many fine, $\gk$-additive
measures on $P_\kappa(\lambda)$.

\end{theorem}

\begin{proof}
As in the proof of Theorem \ref{Theorem.Supercompact}, we
may assume without loss of generality that there are at
least $\gl^+$ many
%carries $2^{2^{[\gl]^{< \gk}}} = 2^{2^\gl} = \gl^{++}$
fine, normal, $\gk$-additive measures on $P_\gk(\gl)$ in
$V$. If we force with Magidor's iteration of Prikry forcing
found in \cite{Ma} which turns $\gk$ into the least
measurable cardinal, then by the work of \cite{Ma}, $\gk$
remains $\gl$ strongly compact. Further,
%(in fact, $\gk$ remains $2^\gl$ strongly compact),
each fine, normal, $\gk$-additive measure ${\cal U}$ on
$P_\gk(\gl)$ extends to a fine, $\gk$-additive measure
${\cal U}^* \supseteq {\cal U}$ on $P_\gk(\gl)$ in the
generic extension $\ov V$ resulting from the iterated
Prikry forcing. Thus, if we now force over $\ov V$ with
$\FP = {\rm Add}(\go, 1) \ast \dot {\rm Coll}(\gl^+,
\gl^{++})$, giving rise to the forcing extension $\ov
V[c][G]$, then the argument given in the proof of Theorem
\ref{Theorem.Supercompact} (using strong compactness
embeddings rather than supercompactness embeddings) goes
through and shows there are $\gl^+$ many fine,
$\gk$-additive measures on $P_\gk(\gl)$ in $\ov V[c][G]$.
%$\ov V^{{\rm Add}(\go, 1) \ast \dot {\rm
%Coll}(\gl^+, \gl^{++})} \models ``P_\gk(\gl)$ carries
Since all partial orderings used preserve the GCH and the
fact $\gk$ is the least measurable cardinal, $\gk$ is $\gl$
strongly compact but isn't $\gl$ supercompact in $\ov
V[c][G]$.
%This completes the proof of Theorem \ref{Theorem.StronglyCompact}.
\end{proof}

We note that since forcing with $\coll(\gl^+, \gl^{++})$ in
both Theorems \ref{Theorem.Supercompact} and
\ref{Theorem.StronglyCompact} adds a new subset of $\gl^+$
to a model obtained by small forcing, it follows by the
main result of \cite[p. 552]{HS} that $\gk$ isn't $\gl^+ =
2^\gl$ strongly compact in $\ov V[c][G]$.
%shows that the forcing $\add(\omega,1) \ast
%\dot \coll(\lambda^+,\lambda^{++})$ definitely
%destroys the $\lambda^+$ strong compactness of $\kappa$.

To this point, in Section \ref{s3}, we have constructed
generic extensions in which we force over a ground model
satisfying the GCH and obtain another universe which also
satisfies the GCH. Our methods, however, also work in
situations where the GCH fails. For instance, we have the
following theorem.

\begin{theorem}\label{NoGCH}
Suppose that $\kappa$ is $\kappa^{++}$ supercompact and the
GCH holds in $V$. Then there is a forcing extension in
which:

\begin{enumerate}\setlength{\parskip}{0pt}

 \item $2^\gd = 2^{\gd^+} = \gd^{++}$ for every inaccessible
       cardinal $\delta\leq\kappa$.

 \item $\kappa$ is $\kappa^+$ supercompact, but $\gk$
         isn't $2^\kappa$ supercompact.

 \item There are exactly $\gk^{++}$ many normal measures on $\gk$.

 \item  There are exactly
       $\gk^{++}$ many fine, normal, $\gk$-additive measures
       on $P_\kappa(\kappa^+)$.

\end{enumerate}

\end{theorem}

\noindent The point of Theorem \ref{NoGCH} is that the GCH
fails at $\gk$, $\gk$ is $\gk^+$ supercompact, yet both
$\gk$ and $P_\gk(\gk^+)$ have fewer than the maximal number
of normal measures.

\medskip
\begin{proof}
To prove Theorem \ref{NoGCH}, suppose we start with a
ground model $V$ in which the GCH holds and $\gk$ is
$\gk^{++}$ supercompact. Force over $V$ with the reverse
Easton iteration of length $\gk + 1$ which begins by adding
a Cohen subset of $\go$ and then does nontrivial forcing
only at those stages
%$\gd \le \gk$
which are inaccessible cardinals in $V$. At such a stage
$\gd \le \gk$, we force with $\add(\gd, \gd^{++})$.
Standard arguments then show that in the resulting model
$\ov V$, for every inaccessible cardinal $\gd \le \gk$,
$2^\gd = 2^{\gd^+} = \gd^{++}$, and $\gk$ is $2^\gk =
2^{{[\gk^+]}^{< \gk}} = 2^{\gk^+} = \gk^{++}$ supercompact.
It is consequently the case that in $\ov V$, there are
exactly $2^{2^\gk} = 2^{\gk^{++}} = \gk^{+++}$ many normal
measures on $\gk$, and there are exactly
$2^{2^{{[\gk^+]}^{< \gk}}} = 2^{\gk^{++}} = \gk^{+++}$ many
fine, normal, $\gk$-additive measures on $P_\gk(\gk^+)$.

If we now force over $\ov V$ with $\FP = \add(\go, 1) \ast
\dot {\rm Coll}(\gk^{++}, \gk^{+++})$, giving rise to the
forcing extension $\ov V[c][G]$, then the proofs of
Theorems \ref{MainTheorem} and \ref{Theorem.Supercompact}
remain valid and show that in $\ov V[c][G]$, there are
exactly $\gk^{++}$ many normal measures on $\gk$, and there
are exactly $\gk^{++}$ many fine, normal, $\gk$-additive
measures on $P_\gk(\gk^+)$. Since forcing with $\FP$
preserves the $\gk^+$ supercompactness of $\gk$, and since
for the same reasons as mentioned after the proof of
Theorem \ref{Theorem.StronglyCompact}, $\gk$ isn't
$\gk^{++}$ strongly compact in $\ov V[c][G]$, the proof of
Theorem \ref{NoGCH} is now complete.
\end{proof}

With a little more work, it is possible to obtain the
conclusions of Theorem \ref{NoGCH} with $\gk$ in addition
being the least measurable cardinal. A brief outline of the
argument is as follows. Start with a ground model $V$
containing a cardinal $\gk$ for which $2^\gk = 2^{\gk^+} =
\gk^{++}$, $2^{\gk^{++}} = \gk^{+++}$, $\gk$ is $\gk^+$
supercompact, and $\gk$ is the least measurable cardinal.
(Note that the consistency of a cardinal with these
properties is originally due to Woodin. A construction of a
model containing such a cardinal may be found, for example,
in \cite{AC1}.) Since $2^\gk = \gk^{++}$ holds in $V$,
there are $\gk^{++}$ many permutations $\pi : \gk \to \gk$
in $V$. Let $j : V \to M$ be an elementary embedding
witnessing the $\gk^+$ supercompactness of $\gk$ such that
$\gk$ isn't measurable in $M$. Force over $V$ with the
reverse Easton iteration $\FP = \FP_\gk \ast \dot \FQ$ of
length $\gk + 1$ which begins by adding a Cohen subset of
$\go$, adds a Cohen subset to each non-measurable
inaccessible cardinal $\gd < \gk$, does trivial forcing at
all other stages $\gd < \gk$, and ends by adding a Cohen
subset of $\gk$. Let $G = G_\gk \ast g$ be $V$-generic over
$\FP$. By a standard argument (see, for example, the proofs
of Lemma 1.1 of \cite{A01} or Lemma 6 of \cite{C93}), in
$V[G_\gk][g]$, there are exactly $2^{2^\gk} = \gk^{+++}$
many fine, normal, $\gk$-additive measures on
$P_\gk(\gk^+)$. Also, for each permutation $\pi : \gk \to
\gk$, a standard argument allows us to lift $j$ to $j^+_\pi
: V[G_\gk][g] \to M[G_\gk][g_\pi][H][g^+]$, where $g_\pi =
\pi''g$. If $i_\pi : V[G_\gk][g] \to N$ is the induced
normal ultrapower embedding, so that $j^+_\pi$ factors as
$k_\pi \circ i_\pi$, then it follows that
$i_\pi(G_\gk)(\gk) = g_\pi$. Since there are $\gk^{++}$
many permutations $\pi : \gk \to \gk$ in $V$, this means
that in $V[G_\gk][g]$, there are (at least) $\gk^{++}$ many
normal measures on $\gk$. If we now force with $\add(\go,
1) \ast \dot {\rm Coll}(\gk^{++}, \gk^{+++})$, the same
argument given in the proof of Theorem \ref{NoGCH} then
shows that in the resulting generic extension, there are
$\gk^{++}$ many normal measures on $\gk$ and $\gk^{++}$
many fine, normal, $\gk$-additive measures on
$P_\gk(\gk^+)$.

The arguments we have presented allow us often to conclude
that a large cardinal has strictly fewer than the maximal
number of measures of the desired kind, even when we are
unable to calculate the exact number of measures. For
instance, if we start with a ground model $V$ where
$\kappa$ is the least measurable cardinal, $\gk$ is $\gk^+$
supercompact, $2^\gk = 2^{\gk^+} = \gk^{++}$, and
$2^{\gk^{++}} = \gk^{+++}$, then after forcing with
$\add(\go, 1) \ast \dot {\rm Coll}(\gk^{++}, \gk^{+++})$,
we may conclude that there are at most $\gk^{++}$ many
normal measures on $\gk$ and at most $\gk^{++}$ many fine,
normal, $\gk$-additive measures on $P_\gk(\gk^+)$. Thus,
even if we don't know either the exact number of normal
measures on $\gk$, or fine, $\gk$-additive measures on
$P_\gk(\gk^+)$, or fine, normal, $\gk$-additive measures on
$P_\gk(\gk^+)$ in our ground model, we still know that
there are fewer than the maximal number of such measures in
the forcing extension.

The method is quite malleable and allows for diverse
similar results. For instance, by starting with a model in
which the GCH holds, $\gk$ is measurable, $\gl > \gk$ is
regular, and for $\gd \le \gk^{++}$ any finite or infinite
cardinal, there are exactly $\gd$ many normal measures on
$\gk$, by forcing with $\add(\go, 1) \ast \dot \add(\gk^+,
\gl)$, we have constructed a model in which $\gk$ is
measurable, $2^\gk = \gk^+$, $2^{\gk^+} = \gl$, and there
are still $\gd$ many normal measures on $\gk$. In addition,
we may force to restrict the number of extenders witnessing
the $\gl$ strongness of $\gk$.

One limitation of the method is that after the forcing
constructions of Theorems \ref{Theorem.Supercompact} and
\ref{Theorem.StronglyCompact}, as we have mentioned, $\gk$
is no longer strongly compact. Thus, it is not possible to
use the techniques of this paper to construct a strongly
compact cardinal $\gk$ such that there are fewer than the
maximal number of fine, $\gk$-additive measures on
$P_\gk(\gl)$, for $\gl > \gk$ regular. In addition, the
arguments of this paper do not seem to allow us to
construct a model in which $\gk$ is the least measurable
cardinal and there are exactly $\gk$ many normal measures
on $\gk$. If one were to use $\add(\omega,1) \ast \dot
\coll(\kappa,2^{2^\kappa})$ in the arguments of Theorem
\ref{MainTheorem} above, or any other small forcing
followed by the collapse of an ordinal to $\kappa$, then
the main theorem of \cite{H6} shows that the measurability
of $\kappa$ would be destroyed.

We conclude with a list of open questions that the methods
of this paper seem not to resolve (although for some, we
now have partial answers). They are as follows:

\begin{enumerate}

\item Is it consistent, relative to anything, for the least
measurable cardinal $\gk$ to have exactly $\gk$ many normal
measures?
%Recall that this question is
%left open by Baldwin in \cite{Ba},
%and as just mentioned, is left
%unresolved by the techniques
%of this paper.

\item How many measures can the least measurable cardinal
have, when there is a strongly compact or supercompact
cardinal above it? Results here show that any regular
cardinal above $\kappa$ is possible; for smaller values, it
seems to be completely open.

\item For which values of $\gl$ does Con(ZFC + There is one
measurable cardinal $\kappa$) imply Con(ZFC + There is a
measurable cardinal $\kappa$ with exactly $\gl$ many normal
measures)? The case $\lambda=1$ is provided by $L[\mu]$ and
the case $\lambda=2^{2^\kappa}$ is provided by the usual
lifting argument techniques. The results of Mitchell
\cite{Mi} and Baldwin \cite{Ba}, however, use measurable
cardinals of high Mitchell order.

\item How many normal measures can a measurable cardinal of
nontrivial Mitchell rank have? In the canonical Mitchell
inner model with $o(\kappa)=\delta$, there are exactly
$|\delta|$ many normal measures on $\kappa$. The usual
lifting arguments show that $2^{2^\kappa}$ is also always
possible, with any value of $o(\kappa)$. Our results here
show that any cardinal $\delta\in [\kappa^+,2^{2^\kappa}]$
with $\cof(\delta)>\kappa$ is also possible with any value
of $o(\kappa)<\kappa^{++}$.

\item How many normal measures can $\gk$ have if $\gk$ is
measurable and $2^\gk > \gk^+$? The work of this paper
shows that such a $\gk$ can have fewer than the maximal
number of normal measures, but does not provide a fully
general answer. Our results show that if $\gk$ is
measurable and $2^\kappa>\kappa^+$, then for any regular
cardinal $\delta$ in the interval
$(\kappa^+,2^{2^\kappa}]$, there is a forcing extension
preserving $2^\kappa>\kappa^+$, where $\kappa$ carries
exactly $\gd$ many normal measures.

\item If $\gk$ is $\gl$ supercompact, how many fine, normal
$\gk$-additive measures can there be on  $P_\gk(\gl)$? If
$\kappa$ is $\lambda$-supercompact, then our results here
show that for any $\delta\in [\lambda^+,2^{2^\lambda}]$
with $\cof(\delta)\geq\lambda^+$, there is a forcing
extension where $\kappa$ has exactly $\delta$ many such
$\lambda$ supercompactness measures.

\item If $\gk$ is $\gl$ strongly compact but isn't $\gl$
supercompact, how many fine, $\gk$-additive measures can
$P_\gk(\gl)$ have? Our results show that any $\delta\in
[\lambda^+,2^{2^\lambda}]$ with $\cof(\delta)\geq\lambda^+$
is possible.

\item For what values of $\gl$ is it consistent for $\gk$
to be fully supercompact and for $\gk$ to have exactly
$\gl$ normal measures not concentrating on measurable
cardinals? The usual lifting arguments show that
$\lambda=2^{2^\kappa}$ is always possible, and so any
$\lambda>\kappa^+$ with $\cof(\lambda)>\kappa^+$ is
possible. For smaller values of $\lambda$, it seems to be
completely open.

\end{enumerate}

\end{document}